\documentclass[11pt,centertags]{amsart}

\usepackage{times}
\usepackage[mathcal]{euscript}

\DeclareSymbolFont{SY}{U}{psy}{m}{n}
\DeclareMathSymbol{\emptyset}{\mathord}{SY}{'306}

\renewcommand{\epsilon}{\varepsilon}

\numberwithin{equation}{section}

\newcommand{\bbC}{\mathbb{C}}

\newcommand{\bbN}{\mathbb{N}}

\newcommand{\bbR}{\mathbb{R}}

\newcommand{\EE}{\mathsf{E}}

\newcommand{\bA}{{\mathbf A}}

\newcommand{\cB}{{\mathcal B}}

\newcommand{\cO}{{\mathcal O}}
\newcommand{\cP}{{\mathcal P}}

\newcommand{\dist}{{\ensuremath{\mathrm{dist}}}}

\newcommand{\fH}{\mathfrak{H}}

\DeclareMathOperator{\Ran}{\mathrm{Ran}}
\DeclareMathOperator{\Ker}{\mathrm{Ker}}

\newcommand{\spec}{{\ensuremath{\mathrm{spec}}}}

\renewcommand{\index}{\mathrm{index}}

\newcommand{\convhull}{\mathrm{conv.hull}}

\newtheorem{theorem}{Theorem}[section]
\newtheorem{lemma}[theorem]{Lemma}
\newtheorem{proposition}[theorem]{Proposition}

\newtheorem{remark}[theorem]{Remark}

\newtheorem{example}[theorem]{Example}{\bf}{\rm}

\newtheorem{introtheorem}{Theorem}{\bf}{\it}
{\bf}{\it}
{\bf}{\it}
\newtheorem{introhypothesis}{Hypothesis}{\bf}{\it}

\title[On a Subspace Perturbation Problem]
{On a Subspace Perturbation Problem}

\author[V. Kostrykin]{Vadim Kostrykin}
\address{Vadim Kostrykin\\ Fraunhofer-Institut f\"{u}r Lasertechnik, Steinbachstra{\ss}e
15, D-52074\\ Aachen, Germany} \email{kostrykin@ilt.fhg.de, kostrykin@t-online.de}

\author[K.A. Makarov]{Konstantin A. Makarov}
\address{Konstantin A. Makarov\\ Department of Mathematics, University of
Missouri, Co\-lum\-bia, MO 65211, USA} \email{makarov@math.missouri.edu}

\author[A.K. Motovilov]{Alexander K. Motovilov}
\address{Alexander K. Motovilov\\ Joint Institute for
Nuclear Research, 141980 Dubna, Moscow Region, Russia} \curraddr{Department of Mathematics,
University of Missouri, Columbia, MO 65211, USA}
\email{motovilv@thsun1.jinr.ru}



\begin{document}

\begin{abstract}
We discuss the problem of perturbation of spectral subspaces for linear self-adjoint
operators on a separable Hilbert space. Let $A$ and $V$ be bounded  self-adjoint
operators. Assume that the spectrum of $A$ consists of two disjoint parts $\sigma$
and $\Sigma$ such that $d=\text{dist}(\sigma, \Sigma)>0$. We show that the norm of
the difference of the spectral projections $\EE_A(\sigma)$ and $\EE_{A+V}\big
(\{\lambda \,\, | \,\,   \dist(\lambda, \sigma)<d/2\}\big)$ for  $A$ and  $A+V$ is less then one
whenever either (i) $\|V\|<\frac{2}{2+\pi}d$ or (ii) $\|V\|<\frac{1}{2}d$ and certain
assumptions on the mutual disposition of the sets $\sigma$ and $\Sigma$ are satisfied.
\end{abstract}

\maketitle

\section{Introduction}

It is well known (see, e.g., \cite{Kato}) that if $A$ and $V$ are bounded
self-adjoint operators on a separable Hilbert space $\fH$, then (the perturbation)
$V$ does not close gaps of length greater than $2\|V\|$ in the spectrum of $A$. More
precisely, if $(a, b)$ is a finite interval and $(a, b)\subset \varrho(A)$, the
resolvent set of $A$, then
\begin{equation*}
(a + \|V\|, b - \|V\|)\subset \varrho(A+sV)
\quad \text{ for all } s\in [-1, 1]
\end{equation*}
whenever $2\|V\|<b-a$. Hence, under the assumption that $A$ has an isolated part $\sigma$ of the
spectrum separated from its remainder by gaps of length greater than or equal to $d>0$, the
spectrum of the operators $A+sV$, $s\in [-1, 1]$ will also have separated components, provided
that the condition
\begin{equation}\label{NGC}
 \|V\|<\frac{d}{2}
\end{equation}
holds.

Our main concern is to study the variation the corresponding spectral subspace associated with the
isolated part $\sigma$ of the spectrum of $A$ under perturbations satisfying \eqref{NGC}.

{}For notational setup we  assume the following hypothesis.
\begin{introhypothesis}\label{h1} Assume that $A$ and $V$ are
  bounded self-adjoint operators
on a separable
Hilbert space $\fH$.
Suppose  that the spectrum of  $A$  has a part $\sigma$ separated from the
remainder of the spectrum
$\Sigma$ in the sense that
\begin{equation*}
 \spec(A)=\sigma\cup\Sigma
\end{equation*}
and
\begin{equation*}
 \dist(\sigma, \Sigma)=d>0.
\end{equation*}
Introduce the orthogonal projections $P=\EE_{A}(\sigma)$ and
$Q=\EE_{A+V}(U_{d/2}(\sigma)), $ where $U_\varepsilon(\sigma)$, $\varepsilon
>0$ is the open $\varepsilon$-neighborhood of the set $\sigma$.
Here $\EE_{A}(\Delta)$ and $\EE_{A+V}(\Delta)$ denote the spectral projections for operators $A$
and $A+V$, respectively, corresponding to a Borel set $\Delta\subset\bbR$\,.
\end{introhypothesis}

In this note we address the following question:
\emph{Assuming  Hypothesis \ref{h1}, does condition \eqref{NGC}
 imply}
\begin{equation*}
\|P-Q\|<1?
\end{equation*}
We give a partially affirmative answer to this question.
 The precise statement reads as follows.

\begin{introtheorem}\label{main} Assume Hypothesis \ref{h1} and suppose that
either
\smallskip

{\rm{(i)}} $ \|V\|<\frac{2}{2+\pi }d$
\smallskip

\noindent or

{\rm{(ii)}} $\|V\|<\frac{1}{2}d$

\noindent and
\begin{equation}
\label{chull}
\convhull(\sigma)\cap\Sigma=\emptyset \,\,\,
\text{\, or \, }\convhull(\Sigma)\cap\sigma=\emptyset.
\end{equation}
Then
\begin{equation*}
 \|P-Q\|<1.
\end{equation*}
\end{introtheorem}

Our strategy of  proof of Theorem \ref{main} does not allow to
relax   condition
\begin{equation}\label{pipi}
\|V\|<\frac{2}{2+\pi}d
\end{equation}
and just assume the natural condition \eqref{NGC} with no
additional hypotheses. It is an
\emph{open problem} whether
 Hypothesis \ref{h1}
alone and the bounds
\begin{equation}\label{pipi2}
 \frac{2}{2+\pi}\le \frac{\|V\|}{d} <\frac{1}{2}
\end{equation}
on the perturbation $V$ imply $\|P-Q\|<1$.

For compact perturbations $V$ satisfying  inequality \eqref{NGC} we can however state that the
pair $(P,Q)$ of the  orthogonal projections is a Fredholm pair with zero index. Recall that the
pair $(P,Q)$ of orthogonal projections is called Fredholm if the operator $QP$ viewed as a map
from $\Ran P$ to $\Ran Q$ is a Fredholm operator \cite{Avron:Seiler:Simon}. The index of this
operator is called the index of the pair $(P,Q).$

\begin{introtheorem}\label{Theorem:2}
Assume Hypothesis \ref{h1} and suppose that $V$ is a compact operator satisfying
\eqref{NGC}. Then the pair $(P,Q)$ is Fredholm with zero index. In particular, the subspaces
$\Ker(PQ^\perp - I)$ and $\Ker(P^\perp Q - I)$ are finite-dimensional and
\begin{equation*}
\dim\Ker(PQ^\perp - I) = \dim \Ker(P^\perp Q - I).
\end{equation*}
\end{introtheorem}

In the ``overcritical" case $\|V\|>d/2$, the perturbed operator $A+V$ may not have separated parts
of the spectrum at all. In this case we give  an example
where the spectral measure of the perturbed operator $A+V$ is ``concentrated'' on the unit sphere
in the space of bounded operators $\cB(\fH)$ centered at the point $P=\EE_A(\sigma)$, with the
norm of the perturbation being arbitrarily close to $d/2$. That is, given $d>0$, for any
$\varepsilon>0$ one can find a self-adjoint operator $A$ satisfying Hypothesis \ref{h1} and a
self-adjoint perturbation $V$ with $\|V\|= d/2+\varepsilon$ such that
\begin{equation*}
\|\EE_A(\sigma)- \EE_{A+V}(\Delta)\|=1
\end{equation*}
for   any Borel set $\Delta\subset \bbR.$

\subsection*{Acknowledgments.} V.~Kostrykin is grateful to V.~Enss
for useful discussions. K.~A.~Ma\-ka\-rov is grateful to F.~Gesztesy for
 critical remarks. A.~K.~Motovilov acknowledges the great hospitality and financial support by the
Department of Mathematics, University of Missouri--Columbia, MO, USA. He was also supported in
part by the Russian Foundation for Basic Research within  the RFBR Project 01-01-00958.

\section{Proof of Theorem \ref{main}}

Our  proof of Theorem \ref{main} is based on the following sharp result (see
\cite{McEachin} and references cited therein) taken from geometric perturbation theory
initiated by C. Davis \cite{Davis:58} and developed further in  \cite{Bhatia:Davis:McIntosh},
\cite{Bhatia:Davis:Koosis}, \cite{Davis:123}, \cite{Davis:Kahan}, \cite{Kato}.

\begin{proposition}\label{mce}
Let $A$ and $B$ be bounded self-adjoint operators and $\delta$ and $\Delta$ two Borel sets on the
real axis $\bbR$\,. Then
\begin{equation*}
\dist (\delta, \Delta)
\| \EE_A(\delta) \EE_B(\Delta)\|
\leq \frac{\pi}{2}
\| A-B \|.
\end{equation*}
If, in addition, the convex hull of the set $\delta$ does not intersect the set
$\Delta$, or the convex hull of the set $\Delta$ does not intersect the set
$\delta$,
then one has the  stronger result
\begin{equation*}
\dist (\delta, \Delta)
\| \EE_A(\delta) \EE_B(\Delta)\|
\leq
\| A-B \|.
\end{equation*}
\end{proposition}

We split the proof of Theorem \ref{main} into the following two lemmas.

\begin{lemma}\label{1}
Assume Hypothesis \ref{h1}. Assume, in addition, that
\eqref{pipi} holds. Then
\begin{equation*}
\|P-Q\|<1.
\end{equation*}
\end{lemma}

\begin{proof}
Clearly
$\spec(A+V)\subset \overline{U_{\|V\|}(\sigma\cup \Sigma)}$, where bar denotes
the (usual) closure in $\bbR$, and then
\begin{equation*}
Q^\perp=\EE_{A+V}\big(\overline{U_{\|V\|}(\Sigma)}\big).
\end{equation*}
By the first claim of Proposition \ref{mce},
\begin{equation}\label{pqort}
\|PQ^\perp\|\le \frac{\pi}{2}\frac{\|V\|}{\dist(\sigma, U_{\|V\|}(\Sigma))}.
\end{equation}
The distance between the set $\sigma$ and the $\|V\|$-neighborhood
of the set
$\Sigma$  can be estimated from below as follows,
\begin{equation*}
\dist(\sigma, U_{\|V\|}(\Sigma))\ge d-\|V\|>0.
\end{equation*}
Then  \eqref{pqort} implies the inequality
\begin{equation*}
\|PQ^\perp\|\le \frac{\pi}{2}\frac{\|V\|}{d-\|V\|}.
\end{equation*}
Hence, from  inequality \eqref{pipi} it follows that
\begin{equation}\label{pq<1}
\|PQ^\perp\|\le \frac{\pi}{2}\frac{\|V\|}{d-\|V\|}<1.
\end{equation}
Interchanging the roles of $\sigma$ and $\Sigma$ one obtains the analogous inequality
\begin{equation}\label{qp<1}
\|P^\perp Q\| < 1.
\end{equation}
Since
\begin{equation}\label{glaz}
\|P-Q\|=\max \{\|PQ^\perp\|, \|P^\perp Q\|\}
\end{equation}
(see, e.g., \cite[Ch. III, Section 39]{Akhiezer:Glazman}),  inequalities
\eqref{pq<1} and \eqref{qp<1} prove the assertion.
\end{proof}

Under additional assumptions on mutual disposition of the parts $\sigma$ and $\Sigma$ of the
spectrum of $A$ one can relax the condition \eqref{pipi} on the norm of perturbation and replace
it by the natural condition \eqref{NGC}.

\begin{lemma}\label{2.3}
Assume Hypothesis \ref{h1} and suppose that condition
\eqref{NGC} holds.
\begin{enumerate}
\item[(i)] If either
$\sigma \cap \convhull(\Sigma)=\emptyset$ or $\convhull(\sigma)\cap
\Sigma=\emptyset$,
then
\begin{equation}\label{asser}
\|P-Q\|<1.
\end{equation}

\item[(ii)] If in addition
the sets $\sigma $ and $\Sigma$  are subordinated, that is,
\begin{equation*}
\convhull(\sigma)\cap
\convhull(\Sigma)=\emptyset,
\end{equation*}
then the following sharp estimate holds
\begin{equation}\label{sharp:est}
\|P-Q\|<\frac{\sqrt{2}}{2}.
\end{equation}
\end{enumerate}
\end{lemma}
\begin{proof} (i) The proof follows that of
Lemma \ref{1}. Applying the second assertion of Proposition \ref{mce} instead of
inequality \eqref{pqort},  one derives the estimates
\begin{equation}\label{pqort2}
\|PQ^\perp\|\le \frac{\|V\|}{\dist(\sigma, U_{\|V\|}(\Sigma))}
\le \frac{\|V\|}{d-\|V\|}<1,
\end{equation}
under hypothesis \eqref{pipi2}, and then the inequality $\|P^\perp Q\|<1$,
  proving  assertion \eqref{asser} using
\eqref{glaz}.

(ii) First assume that $V$ is off-diagonal, that is,
\begin{equation*}
\EE_A(\sigma)V\EE_A(\sigma)= \EE_A(\sigma)^\perp V\EE_A(\sigma)^\perp = 0.
\end{equation*}
Then
the inequality $\|P-Q\|<\frac{\sqrt{2}}{2} $ follows from the $\tan 2 \Theta$-Theorem proven first
by C.~Davis (see, e.g., \cite{Davis:Kahan})
\begin{equation*}
\|P-Q\|\le \sin \left(\frac{1}{2}\arctan
\frac{2\|V\|}{d}\right)<\frac{\sqrt{2}}{2}.
\end{equation*}
A related result can be found in \cite{Adamyan:Langer:95}.

The general case can be reduced to  the off-diagonal one by the following trick.
Assume that $V$
is not necessarily off-diagonal. Decomposing the perturbation $V$ into the diagonal
$V_{\mathrm{diag}}$ and off-diagonal $V_{\mathrm{off}}$ parts with respect to the orthogonal
decomposition  $\fH=\Ran \EE_{A}(\sigma)\oplus \Ran \EE_{\bA}(\sigma)^\perp$
associated with the range of the projection $\EE_A(\sigma)$
\begin{equation*}
V=V_{\mathrm{diag}}+V_{\mathrm{off}},
\end{equation*}
one concludes that
\begin{equation*}
\EE_{A+V_{\mathrm{diag}}}(U_{d/2}(\sigma))=\EE_{A}(\sigma ).
\end{equation*}
Moreover, the distance between
the spectrum of the part of  $A+V_{\mathrm{diag}}$
associated with the invariant subspace
$\Ran \EE_{A+V_{\mathrm{diag}}}(U_{d/2}(\sigma))$ and
the remainder of the spectrum of $A+V_{\mathrm{diag}}$ does not exceed
$d-2\|V_{\text{diag}}\| >
0$. Using the $\tan 2
\Theta$-Theorem then yields
\begin{align*}
 \|P-Q\|&\le \sin
\left(\frac{1}{2}\arctan
\frac{2\|V_{\text{off}}\|}{d-2\|V_{\text{diag}}\|}\right)\\
&\le \sin
\left(\frac{1}{2}\arctan
\frac{2\|V\|}{d-2\|V\|}\right)<\frac{\sqrt{2}}{2},
\end{align*}
completing the proof.
\end{proof}

The sharpness of  estimate \eqref{sharp:est} is shown by the following example.

\begin{example}
Let $\fH=\bbC^2$. For an arbitrary $\epsilon\in(0,3/4)$ consider the  $2\times 2$ matrices
\begin{equation*}
A=\begin{pmatrix} 0 & 0 \\ 0 & 1\end{pmatrix}, \qquad V=\begin{pmatrix} 1/2-\epsilon &
\sqrt{\epsilon}/2 \\ \sqrt{\epsilon}/2 & -1/2+\epsilon\end{pmatrix}.
\end{equation*}
Let $\sigma=\{0\}$ and  $\Sigma=\{1\}$. Obviously, $\dist(\sigma, \Sigma)=1$.
Since
\begin{equation*}
\|V\| = \frac{1}{2}\sqrt{1-3 \epsilon +4 \epsilon^2} <\frac{1}{2},
\end{equation*}
the perturbation $V$ satisfies the hypotheses of Lemma \ref{2.3}.
 Simple calculations
yield
\begin{equation*}
\begin{split}
Q & =\EE_{A+V}\bigl(U_{1/2}(\sigma)\bigr)
=\EE_{A+V}\bigl((-1/2,1/2)\bigr)\\ & =
\frac{1}{1+(2\sqrt{\epsilon}+\sqrt{1+4\epsilon})^2}\begin{pmatrix}
(2\sqrt{\epsilon}+\sqrt{1+4\epsilon})^2 & -2\sqrt{\epsilon}-\sqrt{1+4\epsilon}
\\
-2\sqrt{\epsilon}-\sqrt{1+4\epsilon} & 1\end{pmatrix},
\end{split}
\end{equation*}
and hence,
\begin{equation*}
\|P-Q\| = \left[1+ (2\sqrt{\epsilon}+\sqrt{1+4\epsilon})^2 \right]^{-1/2} < \frac{\sqrt{2}}{2}.
\end{equation*}
Taking $\epsilon$ sufficiently small, the norm $\|P-Q\|$ can be made arbitrarily close to
$\sqrt{2}/2$.
\end{example}

\section{Proof of Theorem \ref{Theorem:2}}

\begin{lemma}\label{Lemma:X}
Assume Hypothesis \ref{h1} and suppose, in addition, that $V$ is a compact
 operator satisfying
condition \eqref{NGC}.
Then there is a unitary $W$ such that $Q=W P W^\ast$ and $W-I$ is compact.
\end{lemma}

\begin{proof}
Fix $\epsilon>0$ such  that $(1+\epsilon)\|V\|<d/2$ and introduce the
family of spectral projections
\begin{equation*}
\cP(s)=\EE_{A+sV}(U_{d/2}(\sigma)),  \quad s\in(-\epsilon, 1+\epsilon).
\end{equation*}
Clearly, $\cP(0)=P$ and $\cP(1)=Q$. From the analytical perturbation theory (see
\cite{Kato}) one concludes that the operator-valued function $\cP(s)$ is real-analytic on
$(-\epsilon,1+\epsilon)$. Moreover (see \cite[Section II.4.2]{Kato}),
\begin{equation*}
\cP(s) = X(s) \cP(0) X(s)^\ast, \quad s\in[0,1],
\end{equation*}
where $X(s)$ is the unique unitary solution to the initial value problem
\begin{align*}
X'(s) & = H(s) X(s), \quad s\in[0,1],\\
X(0)&=I,
\end{align*}
with $H(s)=\cP'(s) \cP(s) - \cP(s) \cP'(s)$.

Let $\Gamma$ be a Jordan counterclockwise oriented contour encircling $U_{d/2}(\sigma)$ in a way
such that no point of $U_{\|V\|}(\Sigma)$ lies within $\Gamma$. Then
\begin{equation*}
\cP(s) = - \frac{1}{2\pi i}
\int_\Gamma (A + s V - z)^{-1} dz,\quad s\in [0,1],
\end{equation*}
and hence,
\begin{equation*}
\cP'(s) = \frac{1}{2\pi i}\int_\Gamma (A + s V - z)^{-1} V (A + s V - z)^{-1} dz,\quad s\in [0,1].
\end{equation*}
By the  hypothesis  $V$ is compact, and hence, $\cP'(s)$,
$s\in [0,1]$ is also compact, which implies that $H(s)$ is a
compact operator for  $s\in[0,1]$.

Applying the successive approximation method
\begin{equation*}
X_n(s) = I + \int_0^s H(t) X_{n-1}(t) dt,\quad X_0(s)=I,
\end{equation*}
yields that  $X_n(s)$ converges to $X(s)$, $s\in [0,1]$ in the norm topology and
$X_n(s)-I$ is compact for all $n\in\bbN$. Thus, $X(s)-I$ is a compact operator for all $s\in
[0,1]$. Taking $W=X(1)$ yields $Q=WPW^*$, completing  the proof.
\end{proof}

Lemma \ref{Lemma:X} implies that the operator $PWP$ viewed as a map from $\Ran P$ to $\Ran P$ is
Fredholm with zero index. By Theorem 5.2 of \cite{Avron:Seiler:Simon} it follows that the pair
$(P,Q)$ is Fredholm and $\index(P,Q)=\index(PW|_{\Ran P})=0$, proving  Theorem \ref{Theorem:2}.

\section{Overcritical perturbations}

If the perturbation $V$ closes a gap between the separated parts $\sigma$ and $\Sigma$ of the
spectrum of the unperturbed operator $A$, then, necessarily, we are dealing with the case $\|V\|
\geq d/2$. In this case one encounters  a new phenomenon: It may happen that any invariant
subspace of the operator $A+V$ contains a nontrivial element orthogonal to $\Ran
P=\Ran\EE_A(\sigma)$.

To illustrate this phenomenon we need the following abstract result.

\begin{lemma}\label{prop:1}
Let $A$ and $V$ be bounded self-adjoint operators and $\sigma\neq\emptyset$ be a
finite set consisting of isolated eigenvalues of $A$ of finite multiplicity.
 Assume
that the spectrum of the operator $A+V$ has no pure point component. Then for the
orthogonal projection $Q$ onto an arbitrary invariant subspace of the operator $A+V$
the subspace $\Ker ( P^\perp Q-I)$, where $P=\EE_A(\sigma)$,
 is infinite-dimensional.
In particular,
\begin{equation}\label{=1}
 \|P-Q\|=1.
\end{equation}
\end{lemma}

\begin{proof}
Since $A+V$ has no eigenvalues, $\Ran Q$ is an
infinite-dimensional subspace. By hypothesis, $\Ran P$ is a
finite-dimensional subspace. Thus, there exists an orthonormal
system $\{f_n\}_{n\in \bbN}$ in  $\Ran Q$ such that $f_n$ is
orthogonal to $\Ran P$ for any $n\in \bbN$ and hence $P^\perp
Qf_n=f_n$, $n\in \bbN$, proving $\dim\bigl(\Ker(P^\perp
Q-I)\bigr)=\infty$. Now equality \eqref{=1} follows from
representation \eqref{glaz}.
\end{proof}

The next lemma  shows that an isolated {\it eigenvalue} of the
unperturbed operator $A$ separated from the remainder of
the spectrum of $A$
by a gap of length $1$ may ``dissolve" in the essential spectrum
of the perturbed operator $A+V$ turning into a
``\emph{resonance}'', with  the norm of the perturbation being
larger but arbitrarily close to $1/2$.

\begin{lemma}
Let $\varepsilon >0$. Let $A$ and $V$ be $2\times 2$ operator matrices in  $\fH=L^2(0,1)\oplus
\bbC$\,,
\begin{equation*}
A=\begin{pmatrix}
M& 0 \\
  0 & -I_\bbC
\end{pmatrix}
\quad\text{and}\quad
V=\begin{pmatrix}
-\left(\frac{1}{2}+\varepsilon\right)I_{L^2(0,1)}  & \sqrt{\varepsilon} v \\
\sqrt{\varepsilon} v^\ast &   (\frac{1}{2}+\varepsilon){I_\bbC}
\end{pmatrix}
\end{equation*}
with respect to the decomposition $\fH=L^2(0,1)\oplus \bbC$\,. Here $M$ denotes the multiplication
operator in $L^2(0,1) $,
\begin{equation*}
(M f)(\mu)= \mu
f(\mu), \quad 0<\mu<1, \quad f\in L^2(0,1),
\end{equation*}
and $v\in \cB\bigl(\bbC, L^2(0,1)\bigr)$
\begin{align}
(vg)(\mu)&=w(\mu)g,\quad \mu\in(0,1), \quad g\in \bbC,\notag \\
w(\mu)&=\sqrt{\mu(1-\mu)}.\notag
\end{align}
If $\varepsilon <2/5$, then the operator $A+V$ has no eigenvalues.
\end{lemma}

\begin{proof}
Assume to  the contrary  that $\lambda \in \bbR$ is an eigenvalue of the perturbed
operator $A+V$, that is,
\begin{equation*}
(\mu -1/2-\varepsilon)f(\mu)+\sqrt{\varepsilon} w(\mu)g=\lambda
f(\mu)\quad \text{a.e. }\mu \in (0,1)
\end{equation*}
and
\begin{equation*}
\sqrt{\varepsilon }\int_0^1d\mu f(\mu)w(\mu) +(-1/2+\varepsilon)g=\lambda
g
\end{equation*}
for some $f\in L^2(0,1)$ and $g\in \bbC$\,. In particular,
\begin{equation*}
f(\mu)=\sqrt{\varepsilon} \frac{w(\mu)}{\lambda-(\mu
-\frac{1}{2}-\varepsilon)}g,
\end{equation*}
and hence $f\notin L^2(0,1)$ whenever $\lambda \in
[-1/2-\varepsilon,1/2-\varepsilon]$ (unless $f=0$ and $g=0$). Thus, the interval
$[-1/2-\varepsilon,1/2-\varepsilon]$ does not intersect the point spectrum of $A+V$.
Moreover, $\lambda\in(-\infty,-1/2-\varepsilon)\cup(1/2-\varepsilon,\infty)$ is an
eigenvalue of $A+V$ if and only if
\begin{equation}\label{eqEx}
\lambda+\frac{1}{2}-\epsilon+\varepsilon \int_0^1 d\mu \frac{\mu(1-\mu)}
{\mu-\frac{1}{2}-\epsilon-\lambda}=0.
\end{equation}
Elementary analysis of the graph of the function on the left-hand side of
 \eqref{eqEx} then yields that under
the condition $0<\varepsilon<2/5$ there is no solution of equation \eqref{eqEx}
in
$(-\infty,-1/2-\varepsilon)\cup(1/2-\varepsilon,\infty)$. Thus, the point spectrum of $A+V$ is
empty.
\end{proof}

\begin{remark}\label{rema}
We note that $\spec(A)=\{-1\}\cup [0,1]$ and hence $\spec(A)$ has two components
separated by a gap of length one, and the norm of the perturbation $V$ may be
arbitrarily close to $1/2$ (from above):
\begin{equation}\label{neploxo}
\|V\|=\sqrt{\left(\frac{1}{\,2\,}+\varepsilon\right)^2+\frac{1}{\,6\,}\,
\varepsilon
}=\frac{1}{2}+\frac{7}{6}\, \varepsilon +\cO(\varepsilon^2) \quad \text{ as }\quad \varepsilon \to 0.
\end{equation}
\end{remark}

Using scaling arguments, Remark \ref{rema} combined with the result of Lemma
\ref{prop:1} shows that given $d>0$, for any $\varepsilon>0$ one can find a
self-adjoint operator $A$ satisfying Hypothesis \ref{h1} and a self-adjoint
perturbation $V$ with $\|V\|= d/2+\varepsilon$ such that
\begin{equation*}
\|\EE_A(\sigma)- Q\|=1
\end{equation*}
for the orthogonal projection $Q$ onto an arbitrary invariant subspace of the
operator $A+V$.


\end{document}